\documentclass[dvips,epsfig,reqno,graphics]{amsart}
\usepackage{amssymb}
\usepackage{amsfonts}
\usepackage{amsmath}
\usepackage{amsxtra}

\title[]
{Certain constant angle surfaces constructed on curves}
\author[]{Ana-Irina Nistor}
\address{
Departement Wiskunde,
Katholieke Universiteit Leuven,
Celestijnenlaan 200B,
3001 Heverlee,
Belgium}
\email{ana.irina.nistor@gmail.com}

\newtheorem{contor}{1.}
\newtheorem{proposition}[contor]{Proposition}
\newtheorem{theorem}[contor]{Theorem}

\def\proof{{\sc Proof.\ }}
\newcommand{\gata}{\hfill\hskip -1cm \rule{.5em}{.5em}}

\def\R{\mathbb{R}}
\def\E{\mathbb{E}}
\def\H{\mathbb{H}}
\def\S{\mathbb{S}}
\def\a{\alpha}

\begin{document}

\maketitle

\input epsfx.tex

\begin{abstract}

In this paper we classify certain special ruled surfaces in $\R^3$ under
the general theorem of characterization of
constant angle surfaces. We study the tangent developable and conical surfaces
from the point of view the constant angle property.
Moreover, the natural extension to normal and binormal constant
angle surfaces is given.

\vspace{1mm}


\vspace{1mm}


\end{abstract}

\section{Introduction}

The developable ruled surfaces represent a special category of ruled
surfaces and are defined as ruled surfaces with vanishing Gaussian
curvature or, moreover, those ruled surfaces which have constant
Gauss map along each ruling (see for details \cite{kn:Kun02}). Since
19th century the developable surfaces captured the attention of
mathematicians and the main properties of these surfaces are
mentioned in almost all monographs and books on Differential
Geometry. Regarding the property of the Gauss map, in \cite{kn:AG04}
is presented an extended study on varieties with degenerate Gauss
maps, meaning that the Gaussian curvature vanishes everywhere, which
includes also the case of developable surfaces.

\smallskip
Due to their flatness isometries with planes are allowed and,
for this reason, recently it became interesting to discover more
ways to use these surfaces in different practical applications. For
example in \cite{kn:RKR} it is proposed a way of constructing and
displaying graphically the kinematic surfaces, a class of surfaces
generated by mutual moving ruled surfaces that touch each other
along their common ruling. A more suggestive application is
described in \cite{kn:GG07}, where some examples of using developable surfaces in
contemporary architecture are discussed.

\smallskip
Motivated again by their flatness property, in this paper we classify the developables and some other
special surfaces constructed on curves from the
point of view of the constancy angle property, i.e. the normal to
the surface makes a constant angle with a fixed direction.
Initially the study of constant angle surfaces was proposed for the product space
$\S^2\times\R$ in \cite{kn:DFVV07}. Then surfaces endowed with this property in other ambient spaces, namely
$\H^2\times\R$, $\R^3$ were investigated (see for details \cite{kn:DM07, kn:DM09, kn:MN08}).

\smallskip
Moreover, in \cite{kn:MN08} it is also given a comparison between the results obtained
so far in the above mentioned ambient spaces.
We mention that in $\S^2\times \R$ one gets positive Gaussian curvature $K$,
in $\H^2\times\R$ negative Gaussian curvature
is obtained while $K$ vanishes identically when the ambient is the Euclidean space.
Explicitly, the constant angle property of the Gauss map with a fixed direction in $\R^3$ is
equivalent to the fact that the Gauss map lies on a circle in the 2-sphere $\S^2$.
As it has no interior points in $\S^2$, the Gaussian curvature of the surface vanishes
identically. It follows that flatness is a consequence of the constancy angle property.
At this point, noticing this common aspect between developable
ruled surfaces and constant angle surfaces in $\R^3$, in the present
paper we would like to see exactly which types of developable
ruled surfaces satisfy the constancy angle property.

It is well known the following classification of surfaces in $\R^3$ involving their degenerate Gauss map: planes,
cylinders, cones and tangent surfaces -- ruled surfaces generated by
the tangent lines in every point of a curve in space. If we drop the flatness property, the construction of
tangent developable surfaces can be used in order to generate some other surfaces constructed on curves
replacing the tangent line to the curve with the normal line or the binormal line respectively,
rising therefore the so called normal surfaces or binormal surfaces.

\smallskip
In next section we mention some basic facts in the general theory of curves and surfaces useful for the
rest of the paper.
In {\em Section $3$} the main result says that {\em the only
tangent developable constant angle surfaces are generated by generalized helices}.
Very recent results implying generalized helices, also called
{\em slope lines}, can be found in \cite{kn:Mun09} (see also its references).

\smallskip
{\em Section $4$} consists of the extension of the study on the constancy angle property for
the normal and binormal surfaces. The main result is that
{\em the normal constant angle surfaces are pieces of planes}
and {\em the binormal constant angle surfaces are pieces of cylinders}. In {\em Section $5$} we make
some observations regarding the {\em conical constant angle surfaces},
motivated by their affiliation to the classification of flat surfaces in $\E^3$.

\section{ Preliminaries}

Traditionally, the differential geometry of curves starts with a smooth map
of $s$, let's call it $\alpha:I\subset\mathbb{R}\rightarrow \mathbb{R}^3$,
that parameterizes a spatial curve denoted again with $\alpha$. We say that the curve is
parameterized by arc length if $|\alpha'(s)|=1$, where $\alpha'$ is the derivative of $\alpha$ w.r.t. $s$.
Throughout this paper $s$ is the arc length parameter. Let us denote $t(s)=\alpha'(s)$
the (unit) tangent to the curve.
By definition, the curvature of $\alpha$ is $\kappa(s)=|\alpha''(s)|$. If $\kappa\neq 0$,
then the (unit) normal of $\alpha$ can be obtained from
$\alpha''(s)=\kappa(s)n(s)$. Moreover, $b(s)=t(s)\times n(s)$ is called
the (unit) binormal to $\alpha$. With these considerations
$t,\ n,\ b$ define an orthonormal basis. Recall the Frenet-Serret formulae:
\begin{equation}
\label{ain:frenet}
\begin{array}{ll}
t'(s)= \kappa(s) n(s)\\
n'(s)= -\kappa(s) t(s)+\tau(s) b(s)\\
b'(s)= -\tau(s) n(s)
\end{array}
\end{equation}
where $\tau(s)$ is the torsion of $\alpha$ at $s$.

\smallskip
A natural extension from curves to the theory of surfaces constructed on curves can be made as follows.
Given a curve $\alpha$ parameterized by arc length in Euclidean 3-space,
we can think of constructing ruled surfaces involving $\alpha$
and the tangent, normal or binormal lines to the curve.
As a consequence, we have three well known types of surfaces of this kind, namely

\begin{itemize}
\item tangent developable surface: $r(s,v)=\alpha(s)+v t(s)$
\item normal surface: $r(s,v)=\alpha(s)+v n(s)$
\item binormal surface: $r(s,v)=\alpha(s)+v b(s)$
\end{itemize}
(see for details \cite{kn:GAS}). The curve $\alpha$ is the generating curve and the rulings are respectively
the tangent, the normal and the binormal lines to the curve.

\smallskip
The characterization of constant angle surfaces in $\R^3$ was given in \cite{kn:MN08},
where the constant angle is denoted by $\theta$ and without loss of generality, the fixed direction is taken to be
the third real axis, denoted by $k$. The main result is the following

\smallskip
{\bf Theorem A} (\cite{kn:MN08}) {\em  A surface $M$ in $\mathbb{R}^3$ is a constant angle surface if and
only if it is locally isometric to one of the following surfaces:

$(i)$ a surface given by
\begin{equation}
\label{ain:eq1}
r:M \rightarrow \mathbb{R}^{3},\ (u_1,u_2)\mapsto (u_1\cos\theta(\cos u_2,\sin u_2)+\gamma(u_2), u_1\sin\theta)
\end{equation}
with
\begin{equation}
\label{ain:eq2}
\gamma(u_2)=\cos\theta\ \left(-\int\limits_0^{u_2} \eta(\tau)\sin\tau d\tau,
\int\limits_0^{u_2}\eta(\tau)\cos\tau d\tau\right)
\end{equation}
for $\eta$ a smooth function on an interval $I\subset \R$,\\[1mm]
$(ii)$ an open part of the plane $x\sin\theta-z\cos\theta=0$,\\[1mm]
$(iii)$ an open part of the cylinder $\beta\times\mathbb{R}$,
where $\beta$ is a smooth curve in $\mathbb{R}^2$.
}

\smallskip
In the following sections we deal with surfaces constructed on curves and we study their constancy angle property.
We also show how one can retrieve this type of surfaces constructed on curves
from the general theorem of characterization above mentioned.

\section{Tangent developable surfaces}

We start this section with some important properties
of developable ruled surfaces. In \cite{kn:Kun02} it is proved that the parametrization of
every flat ruled surface generically written $(u,v)\mapsto r(u,v)$, in other words an open and dense subset of every flat ruled surface
can be subdivided into
subintervals such that the parametrization corresponding to these subintervals can be included
in one of the following types: plane, cylinder, cone, tangent developable surface.

\smallskip
Thinking now of the class of constant angle surfaces in $\R^3$,
which are flat as we have mentioned in {\em Preliminaries}, we show how one can retrieve the case
of tangent developable surfaces, that is not mentioned
explicitly in {\bf Theorem A} from \cite{kn:MN08},
among the constant angle surfaces in $\mathbb{R}^3$.

\smallskip

First we state and prove the following result concerning which types of tangent
developable surfaces satisfy the constancy angle property:

\begin{theorem}
The tangent developable constant angle surfaces are generated by
cylindrical helices.
\end{theorem}

\proof
Let us consider a tangent developable surface $M$ oriented, immersed in $\mathbb{R}^3$ given by
\begin{equation}
\label{ain:eq3}
r(s,v)=\alpha(s)+v t(s),
\end{equation}

where $\alpha:I\subset\mathbb{R}\rightarrow\mathbb{R}^3$ is a spatial curve
parameterized by arc length consisting of the edge of regression of $M$ and $t$ is the unit tangent to $\alpha$.
The surface $M$ is smooth everywhere, except in points of the curve $\alpha$.

Let us determine the normal to the surface. To do this, we compute the partial derivatives of $r$ with respect to $u$ and $v$
\begin{center}
$r_s(s,v)=\alpha'+v\alpha''\ $ and $\ r_v(s,v)=\alpha'$.
\end{center}

Using now \eqref{ain:frenet}, the normal to the surface is given by
$$
N=\pm \frac{r_s\times r_v}{|r_s\times r_v|}=\mp\ b.
$$
Choosing an orientation of the surface we take
the normal to the surface equal to the binormal of the generating curve $\alpha$.
In the case of constant angle surfaces it follows that the binormal $b$ of
the curve $\alpha$ makes a constant angle $\theta$ with the fixed direction $k$, namely
\begin{equation}
\label{ain:cond}
\widehat{(b,k)}=\widehat{(N,k)}=\theta,\ \theta\in[0,\pi).
\end{equation}

It follows that $\alpha$ is a cylindrical helix.

\gata

\bigskip

In order to write the parametrization of the cylindrical helix $\alpha$ we proceed this way.
In general, for a curve $\alpha$ with $|\alpha'(s)|=1$ and satisfying \eqref{ain:cond} one can write
\begin{equation}
\label{ain:hel}
\displaystyle \alpha(s)=\left(\psi(s),\frac{\mathfrak{b}}{\mathfrak{c}}\ s \right)
\end{equation}
where the curve $\psi:I \subset \R \rightarrow \R^2$
satisfies $\displaystyle |\psi'(s)|=\frac{\mathfrak{a}}{\mathfrak{c}}$ such that
$\mathfrak{a}^2+\mathfrak{b}^2=\mathfrak{c}^2$, $\mathfrak{a,b,c} \in (0,\infty)$.
It results that the derivative of $\psi$ w.r.t. $s$ can be written

$\displaystyle  \psi'(s)=\left(\mathfrak{\frac{a}{c}}\sin\lambda(s), \mathfrak{\frac{a}{c}}\cos\lambda(s) \right)$ for a certain
function $\lambda$.
Integrating, the expression of $\psi$ is obtained, and substituting it
in \eqref{ain:hel} one gets:

\bigskip

The curve $\alpha$ is said to be a {\em cylindrical helix} if it can be parameterized by
\begin{equation}
\label{ain:eq9}
\alpha(s)=\left(\mathfrak{\frac{a}{c}}\int\sin\lambda(s)ds,\ \mathfrak{\frac{a}{c}}\int\cos\lambda(s)ds,\ \mathfrak{\frac{b}{c}}\ s\right),
\end{equation}
where $\mathfrak{a,b,c}$ satisfy the above condition and $\lambda:I\subset\R\rightarrow \R$ is a smooth function.

\smallskip

We would like to see the direct connection between the
tangent developable surfaces  satisfying the constant angle property and {\bf Theorem A}. We have

\begin{theorem}
The tangent developable constant angle surfaces are obtained for
$\eta(\tau)=-\lambda^{-1}(\frac{\pi}{2}-\tau)$ in {\bf Theorem A}.
\end{theorem}

\proof 
We will start the proof with the classical case of a circular helix, namely for $\lambda(s)=-s$, obtaining the parametrization:
\begin{equation}
\label{ain:eq4}
\alpha(s)=\left(\mathfrak{\frac{a}{c}}\cos s, \mathfrak{\frac{a}{c}} \sin s, \mathfrak{\frac{b}{c}} s \right).
\end{equation}

Substituting $\displaystyle \alpha'(s)=\left(-\mathfrak{\frac{a}{c}}\sin s, \mathfrak{\frac{a}{c}} \cos s, \mathfrak{\frac{b}{c}} \right)$
and expression \eqref{ain:eq4} in parametrization \eqref{ain:eq3} we get that the tangent developable of the cylindrical helix $\a$ has the form

\begin{equation}
\label{ain:eq5}
r(s,v)=\left(\mathfrak{\frac{a}{c}}(\cos s - v \sin s), \mathfrak{\frac{a}{c}}( \sin s+ v \cos s), \mathfrak{\frac{b}{c}}(s+v) \right).
\end{equation}

\medskip

We prove that this parametrization is a particular case of item $(i)$ in {\bf Theorem A}.
We determine the general function $\eta$ starting with parametrization \eqref{ain:eq5}
and rewriting it in the form \eqref{ain:eq1}.

\medskip
First, we look at the third component of the parameterizations \eqref{ain:eq1} and \eqref{ain:eq5}.
Recall that the fixed direction $k$ can be decomposed into
its normal and tangent parts and developing the same technique as in
\cite{kn:MN08}, we get that
$$
k=\sin\theta \alpha'+\cos\theta N.
$$

Computing $<r_s,k>$ in two ways, first $<r_s,k>=<\alpha', \sin\theta\alpha'>=\sin\theta$ and secondly
$\displaystyle <r_s,k>=\mathfrak{\frac{b}{c}}$, where $<~,~>$ denotes the Euclidean scalar product,
one gets $\displaystyle \mathfrak{\frac{b}{c}}=\sin\theta$. We
obtain also $\displaystyle\mathfrak{\frac{a}{c}}=\cos\theta$.

\smallskip

After the change of parameter $u_1:=s+v$ in \eqref{ain:eq5}, we get the equivalent parametrization
\begin{eqnarray}
\label{ain:eq6}
r(s,v)  & = & \left(\cos\theta(\cos s -u_1 \sin s + s \sin s)\right.,  \\
   &  & \hspace{1mm}\left. \cos \theta (\sin s + u_1 \cos s - s \cos s), u_1 \sin\theta \right). \nonumber
\end{eqnarray}

A second reparametrization, namely $u_2:= s+\frac{\pi}{2}$, yields

\begin{equation}
\label{ain:eq7}
r(u_1,u_2)=(u_1\cos\theta(\cos u_2,\sin u_2)+\gamma(u_2), u_1 \sin \theta),
\end{equation}
where
\begin{equation}
\label{ain:eq8}
\gamma(u_2)=\cos\theta\left(\sin u_2-\left(u_2-\frac{\pi}{2}\right)\cos u_2, -\cos u_2 - \left(u_2-\frac{\pi}{2}\right)\sin u_2\right).
\end{equation}

Now, by comparison with {\bf Theorem A}, \eqref{ain:eq7} is identically with \eqref{ain:eq1}
and we only have to determine the smooth function $\eta$ in order to write \eqref{ain:eq8} in the same manner with \eqref{ain:eq2}.

\medskip
We claim that $\eta(\tau)=\frac{\pi}{2}-\tau$.

In order to prove the claim we compute
$$
-\int\limits_{0}\limits^{u_2}\eta(\tau)\sin\tau d\tau=
\left(\frac{\pi}{2}-u_2\right)\cos u_2 + \sin u_2-\frac{\pi}{2}
$$

and
$$
\int\limits_{0}\limits^{u_2}\eta(\tau)\cos\tau d\tau=
\left(\frac{\pi}{2}-u_2\right)\sin u_2 - \cos u_2+1.
$$

We complete the proof in this case concluding that the expression \eqref{ain:eq8} is equivalent with \eqref{ain:eq2} for
$\eta(\tau)=\frac{\pi}{2}-\tau$ and taking into account the integration limits, i.e. a translation in the $xOy$ plane.

\medskip

Let us return to the general case of a cylindrical helix.

\smallskip
Follow-on the same idea like in the previous case, already knowing that
$\displaystyle \mathfrak{\frac{a}{c}}=\cos\theta$,
$\displaystyle \mathfrak{\frac{b}{c}}=\sin\theta$ and taking the derivative of $\alpha$ with respect to $s$,

$
\displaystyle \alpha'(s)=\left(\mathfrak{\frac{a}{c}}\sin \lambda(s),\ \mathfrak{\frac{a}{c}}\cos\lambda(s),\ \mathfrak{\frac{b}{c}}\right),
$
we get the parametrization for the tangent developable corresponding to the generalized helix
given by \eqref{ain:eq9} in the form

\begin{eqnarray}
\label{ain:eq10}
r(s,v)  & = & \left(\cos\theta \left(\int \sin \lambda(s)ds+ v \sin\lambda(s)\right)\right.,  \\
   &  & \hspace{1mm}\left. \cos\theta \left(\int \cos \lambda(s)ds+ v \cos\lambda(s)\right), (s+v) \sin\theta \right). \nonumber
\end{eqnarray}

Making now the change of parameters
$u_1:= s+v$ and $u_2:=\frac{\pi}{2}-\lambda(s)$ in \eqref{ain:eq10},
one gets the equivalent parametrization
\begin{equation}
\label{ain:eq11}
r(u_1,u_2)=(u_1\cos\theta(\cos u_2,\sin u_2)+\gamma(u_2), u_1\sin\theta)
\end{equation}
where
\begin{eqnarray}
\label{ain:eq12}
\gamma(u_2)  & = & \cos\theta\left(\int
\left(\lambda^{-1}\big(\frac{\pi}{2}-u_2\big)\right)'\cos u_2\ du_2-
\lambda^{-1}\big(\frac{\pi}{2}-u_2\big)\sin\big(\frac{\pi}{2}-u_2\big)\right., \nonumber \\
  &  & \left. \int
\left(\lambda^{-1}\big(\frac{\pi}{2}-u_2\big)\right)'\sin u_2\ du_2-
\lambda^{-1}\big(\frac{\pi}{2}-u_2\big)\cos\big(\frac{\pi}{2}-u_2\big)\right).
\end{eqnarray}

Our attention is focused on the function $\eta$ which should be determined.
By straightforward computations,

\begin{eqnarray}
-\int\limits_{0}\limits^{u_2}\eta(\tau)\sin\tau d\tau & = & \int_0^{u_2}\eta(\tau)\left(\sin\big(\frac{\pi}{2}-\tau\big)\right)'d\tau\\ \nonumber
& = & \eta(u_2)\sin\big(\frac{\pi}{2}-u_2\big)-\eta(0)-\int\limits_0\limits^{u_2}\eta(\tau)'\cos\tau d\tau .\nonumber
\end{eqnarray}
and
\begin{eqnarray}
\int\limits_{0}\limits^{u_2}\eta(\tau)\cos\tau d\tau & = & - \int_0^{u_2}\eta(\tau)\left(\cos\big(\frac{\pi}{2}-\tau\big)\right)'d\tau\\ \nonumber
& = & - \eta(u_2)\cos\big(\frac{\pi}{2}-u_2\big)+\int\limits_0\limits^{u_2}\eta(\tau)'\sin\tau d\tau .\nonumber
\end{eqnarray}

Taking into account the integration limits
we conclude the proof of the theorem with the fact that also in this general case for $\eta(\tau)= -\lambda^{-1}\big(\frac{\pi}{2}-\tau\big)$
the expressions \eqref{ain:eq12} and \eqref{ain:eq2} are equivalent.

\gata

\bigskip

\section{Normal and Binormal Surfaces}

In this section we deal with the other two types of surfaces constructed
on a spatial curve $\alpha$, {\em normal} and {\em binormal surfaces}.
As we have seen in {\em Preliminaries}, these surfaces are constructed
by using the same technique as the tangent developable ones, replacing in \eqref{ain:eq3} the tangent line $t$
with the normal line $n$, respectively the binormal line $b$.

\medskip
In the first part of this section we study under which conditions
the normal and the binormal surfaces can be retrieved form
{\bf Theorem A} under the property of constant angle surfaces.

Concerning this aspect, we state and prove the following result

\begin{theorem} $ $
\begin{enumerate}
\item The normal constant angle surfaces are pieces of planes.
\item The binormal constant angle surfaces are pieces of cylindrical surfaces.
\end{enumerate}
\end{theorem}

\proof $ $


Let us consider first the parametrization of a normal surface
\begin{equation}
\label{ain:eq14}
r(s,v)=\alpha(s)+v n(s).
\end{equation}

Computing the normal to the above surface, one gets

\begin{equation}
\label{ain:nn}
N=\frac{(1-\kappa v)b-\tau v t}{\sqrt{\Delta}},\ {\rm where}\ \Delta=(1-\kappa v)^2+\tau^2v^2.
\end{equation}

We are interested in those normal surfaces for which the normal
$N$ makes a constant angle $\theta$ with the
fixed direction $k$, namely $(\widehat{N,k})=\theta$, equivalently, $<N,k>=\cos\theta$.
Substituting \eqref{ain:nn} in this expression we get a vanishing polynomial expression of second order in $v$.
So, all the coefficients must be identically zero, that is the following relations are satisfied:
\begin{eqnarray}
&  & \label{ain:n1}<b, k>^2-\cos^2\theta=0 \\
&  & \label{ain:n2}\kappa<b, k>^2+\tau<b, k><t, k>-\kappa\cos^2\theta=0 \\
&  & \label{ain:n3}(\kappa<b, k>+
\tau<t, k>)^2-(\kappa^2+\tau^2)\cos^2\theta=0
\end{eqnarray}

From \eqref{ain:n1} results
\begin{equation}
\label{ain:n11}
<b, k>=\pm\cos\theta.
\end{equation}

Substituting \eqref{ain:n1} in \eqref{ain:n2} we get
 $$\tau<b, k><t, k>=0.$$

Let us distinguish the following cases:

\medskip
\begin{itemize}
\item[a)] $\tau=0.$

\smallskip
Both \eqref{ain:n2} and \eqref{ain:n3} reduces to \eqref{ain:n1} which is automatically fulfilled because $\alpha$ being a plane curve
it's binormal coincides with the normal of the plane. Thinking now the
normal surface like a ruled surface for which the rulings are the normal lines to
the generating plane curve $\alpha$, we get that the
normal constant angle surface is a portion of plane.
\medskip

\item[b)] $\tau\neq 0.$

\smallskip
\begin{itemize}
\item[b.1)] $<b, k>=0.$
From \eqref{ain:n1} follows $\cos\theta=0$ and substituting this result in \eqref{ain:n3} it results $<t, k>=0.$ More,
taking the derivative with respect to $s$ we get also $<n, k>=0.$ We have a contradiction: $k$ is orthogonal to all
$t,\ n,\ b$ which already define an orthonormal basis!

\smallskip
\item[b.2)] $<t, k>=0.$
Analogously we get a contradiction similar to the previous subcase. Again this situation cannot occur.
\end{itemize}
\end{itemize}

\medskip
So, in the case of normal constant angle surfaces we retrieve here the case $(ii)$ from {\bf Theorem A}.


\medskip

In the same manner, consider the parametrization of a binormal surface
\begin{equation}
\label{ain:eq14}
r(s,v)=\alpha(s)+v b(s).
\end{equation}

The normal to the surface is

$$
N=\frac{-n-\tau v t }{\sqrt{\Delta}},\ {\rm where}\ \Delta=1+\tau^2 v^2.
$$

With the same technique as before, one gets that the only case which
can occur is again for $\tau=0$ with the
additional condition $<t, k>=0$. Yet, the binormal to the plane
curve $\alpha$ is parallel to the fixed direction and
$\theta=\frac{\pi}{2}$. Thus, the binormal constant angle surfaces are cylindrical
surfaces generated by the plane curve $\alpha$.
In this manner we reached item $(iii)$ of {\bf Theorem A}.

\gata

\medskip

We conclude this section pointing out that studying these types of surfaces constructed on curves under the constancy angle
property, all three cases in {\bf Theorem A} are retrieved.
First item $(i)$ has as a particular case the tangent developable constant angle surfaces,
the second item $(ii)$ includes the case of normal
surfaces end the last one, $(iii)$ correspond also for binormal surfaces satisfying the constancy angle property.

\bigskip

\section{Conical constant angle surfaces}

In this last section let us return to the classification of flat surfaces in $\R^3$.
In previous sections we recovered the planes, the cylinders and the tangent
developable surfaces among the constant angle surfaces. Considering now the case of conical surfaces regarded
from de point of view of constant angle surfaces, we state and prove

\begin{proposition}
The only conical constant angle surfaces are circular cones.
\end{proposition}

\proof
A conical surface with the vertex in the origin is given by
\begin{equation}
\label{ain:eq31}
r(s,v)=v\alpha(s),
\end{equation}
where we consider now $s,v$ standard parameters. This means that any cone
generated by a generic curve $\alpha$ can be reparameterized using standard
parameters such that $\alpha$ lyes on the unit 2-sphere, $|\alpha(s)|=1$. In these conditions, the normal to the surface is given by
$\displaystyle N=\alpha\times\alpha'$. Constant angle property \eqref{ain:cond} is equivalent with $<\alpha\times\alpha',k>=\theta$.
Deriving with respect to $s$ we get
\begin{equation}
\label{ain:eq32}
<\alpha\times\alpha'', k>=0.
\end{equation}
Now, $\alpha''$ can be decomposed in the orthonormal basis $\{\alpha, \alpha',\alpha\times\alpha'\}$ as
$$\alpha''=<\alpha'', \alpha>\alpha+<\alpha'', \alpha'>\alpha'+<\alpha'', \alpha\times\alpha'>\alpha\times\alpha'.$$

From $|\alpha(s)|=1$ and
$|\alpha'(s)|=1$ we get $<\alpha'',\alpha'>=0$ and $<\alpha'',\alpha>=-1$. Substituting these expressions in the decomposition of $\alpha$
we get
\begin{equation}
\label{ain:eq33}
\alpha''=-\alpha+\kappa_g(\alpha\times\alpha'),
\end{equation}
 where $\kappa_g$ denotes the geodesic curvature of $\alpha$.

Substituting \eqref{ain:eq33} in \eqref{ain:eq32} one obtains that a conical surface is constant angle surface if
verifies
$$\kappa_g<\alpha'\times k>=0.$$ At this point we can conclude that $\alpha$ is a planar curve.
Moreover, knowing that $\alpha$ is on the unit 2-sphere,
it follows that $\alpha$ is a circle. So, \eqref{ain:eq31} parameterizes a circular cone.

\gata

\medskip
In order to retrieve circular cones from {\bf Theorem A}, it suffices to consider a general parametrization for a circular cone and after
a change of frame such that the cone's axis is parallel with the $Oz$-axis, we get the corresponding parametrization from
\eqref{ain:eq2} for $\eta(\tau)=0$ in \eqref{ain:eq3}.

\bigskip

{\small {\bf Acknowledgement.} The author wishes to express her gratitude
to Prof. F. Dillen, {\em Katholieke Universiteit Leuven}, Belgium
 and Prof. M.I. Munteanu, University {\em "Alexandru Ioan Cuza"} of Iasi, Romania, for suggesting her the
subject and for helpful hints given during the preparation of this work.}

\end{document}